\documentclass{rsl}

\gridframe{N}
\usepackage{amsmath,amssymb,bm,natbib,url}
\usepackage{isabelle,isabellesym}

\volume{0}
\issue{0}

\pyear{2014}
\pmonth{Month}
\doinu{}
\title[Review of Symbolic Logic]
      {A Machine-Assisted Proof of G\"odel's Incompleteness Theorems
       for the Theory of Hereditarily Finite Sets}

\author[Lawrence C. Paulson]{LAWRENCE C. PAULSON}
\affil{University of Cambridge}

\leftrunninghead{Lawrence C. Paulson}
\rightrunninghead{G\"odel's Incompleteness Theorems}

\newcommand\SUCC{\mathop{\rm SUCC}}
\newcommand\QR{\mathop{\rm QR}}
\newcommand\Pf{\mathop{\rm Pf\,}}
\newcommand\supp{\mathop{\rm supp}}
\newcommand\quot[1]{\ulcorner#1\urcorner}
\newcommand\tuple[1]{\langle#1\rangle}
\newcommand\vquot[1]{\lfloor#1\rfloor}
\newcommand\ex[2]{\exists#1\,[#2]}
\newcommand\all[2]{\forall#1\,[#2]}
\let\bimp=\leftrightarrow
\let\ts=\thinspace
\renewcommand{\isasymlceil}{\isamath{\ulcorner}}
\renewcommand{\isasymrceil}{\isamath{\urcorner}}

\begin{document}
\maketitle

\begin{abstract}
A formalisation of G\"odel's incompleteness theorems using the Isabelle proof assistant
is described. This is apparently the first mechanical verification of the second incompleteness theorem.
The work closely follows \cite{swierczkowski-finite}, who gave a detailed proof
using hereditarily finite set theory. 
The adoption of this theory is generally beneficial, but it poses certain technical issues
that do not arise for Peano arithmetic.
The formalisation itself should be useful to logicians,
particularly concerning the second incompleteness theorem, where existing proofs are lacking in detail.
\end{abstract}

\section{Introduction.}
G\"odel's incompleteness theorems \citep{goedel-I,goedel31} are undoubtedly the most misunderstood results in mathematics.
\cite{franzen-guide} has written an entire book on this phenomenon. One reason is they have attracted the
attention of a great many non-mathematicians, but even specialists who should know better have drawn
unfounded conclusions. One of the main obstacles to understanding these theorems is the great technical
complexity of their proofs, and indeed of their very statements.

\cite{swierczkowski-finite} claims that the theory of hereditarily finite sets (HF) is more suitable than
the usual Peano Arithmetic (PA) as a basis for proving the incompleteness theorems. The coding of terms and
formulas can be done directly using traditional set-theoretic constructions, without referring to prime
factorisation or the Chinese remainder theorem. As evidence, he gives a detailed presentation of the proofs
of these theorems, along with a development of the HF theory itself. 
He also states a theorem saying that the theories HF and PA are definitionally equivalent.

The present paper describes a
formalisation of {\'S}wierczkowski's development using the interactive theorem prover Isabelle/HOL\@. This
formalisation makes some of the advantages and drawbacks of his approach very clear, and these will be
discussed below. Moreover, the availability of this formal proof (which can be surveyed by anybody
who has a suitable computer and a copy of the Isabelle software) can help to demystify the incompleteness
theorems.

\cite{boolos-provability} devotes more than two pages (pp.\ts33--34) to an explanation of how coding syntax
using integers differs from using PA to reason about addition and multiplication. As a computer scientist, I
do not see the need for such lengthy explanations: coding one thing in another is how computers
work on every architectural level. Coming from that perspective, it isn't obvious that representing the
ordered pair $\tuple{x,y}$ set-theoretically as $\{\{x\},\{x,y\}\}$ is more natural than representing
it arithmetically as $2^x 3^y$, for example. What we can objectively say is that the former approach is
likely to save effort, eliminating the need to formalise the fundamental theorem of arithmetic or
the Chinese remainder theorem explicitly in PA\@.

\clearpage 

It's clear that G\"odel regarded the need to construct explicit formal proofs as highly undesirable. We can
regard the proof of a sentence $A$ on three levels: informally, as a proof of ${}\vdash A$ in a
suitable formal calculus, or as a proof of ${}\vdash \Pf {\quot {A}}$, given a suitable coding system
defining $\quot {A}$ and a provability predicate $\Pf$ corresponding to the formal calculus and coding
system. Obviously, the effort required to prove~$A$ increases hugely as we move up from one level to the
next, but one could argue that the intrinsic complexity does not increase at all; the additional effort is
essentially mechanical and bureaucratic. Nevertheless, G\"odel's treatment makes strenuous efforts to
minimise the need to construct formal proofs.

G\"odel describes a relation $R(x_1,\ldots,x_n)$ as \emph{entscheidungsdefinit} (the modern term is
\emph{numeralwise expressible}) provided there is
a formula ${\rm R}(x_1,\ldots,x_n)$ such that, for each $x_1$, \ldots, $x_n$,
\begin{align}
R(x_1,\ldots,x_n) &\quad\text{implies}\quad {}\vdash {\rm R}(\bm{x_1,\ldots,x_n})
\label{eqn:true} \\
\overline R(x_1,\ldots,x_n) &\quad\text{implies}\quad {}\vdash \neg{\rm R}(\bm{x_1,\ldots,x_n})
\label{eqn:false}
\end{align}
Here, $\overline R$ means ``not $R$" and $\bm{x_1,\ldots,x_n}$ denotes the numerals expressing the values of
$x_1$, \ldots, $x_n$ \cite[p.\ts130]{goedel-I}. This technique shows that ${}\vdash {\rm R}
(\bm{x_1,\ldots,x_n})$ is a theorem of the formal calculus without requiring an explicit proof.
Unfortunately, the price is a considerable increase in intrinsic complexity: explicit numerical bounds have
to be given for all quantifiers, and the proofs that these bounds are sufficiently large can be very
complicated. These proofs refer to the coding functions and require detailed reasoning about primes, lowest
common multiples, etc.

A $\Sigma_1$ formula in PA is logically equivalent to one of the form 
$\exists x_1 \ldots \exists x_n \phi$, where $\phi$ is a primitive recursive formula.
Based on this concept (henceforth simply ``$\Sigma$ formulas''), one can eliminate the need 
for bounded existential quantifiers. 
$\Sigma$ formulas turn out to be sufficient to express the provability predicate $\Pf$ and
the syntactic concepts underlying it: terms, formulas, substitutions, etc. They satisfy property
(\ref{eqn:true}) above but not (\ref{eqn:false}). To recover the latter property, 
\cite{boolos-provability} uses the concept of a $\Delta$ formula:
a $\Sigma$ formula whose negation is also a $\Sigma$ formula. Unfortunately, this approach
again requires bounds for existential quantifiers. \cite{boolos-provability}
devotes more than a page (page~41) to a ``grisly'' proof of one of these bounds, concerned with the coding
of terms.  The very statement of the theorem (which replaces one unbounded existential quantifier by three
bounded quantifiers) is highly technical. As there are a great many other existential quantifiers in the
definition of the provability predicate, this approach cannot lead to an intelligible proof of the
incompleteness theorems.

\cite{swierczkowski-finite} confines himself to $\Sigma$ formulas. Since property (\ref{eqn:false}) does not
hold, it is necessary to perform some proofs in the HF formal calculus. He presents detailed proofs that the
coded substitution operations on coded terms and formulas are single-valued. These proofs are as long as the one
given in \cite{boolos-provability}, but conceptually they are simple; their purpose is to 
demonstrate that the proof of the single-valued property is elementary enough to be proved 
in the HF calculus.

To actually exhibit a formal proof, some elementary concepts and lemmas in the theory of HF have to
be developed formally: the principle of mathematical induction, the linear ordering for the natural numbers,
etc. But to reach the first incompleteness theorem, these formal developments do not even need to define
addition. To reach the second theorem, we require a few addition laws and some basic properties of finite
sequences, but nothing more: certainly, not multiplication. This is the main benefit of using HF, since
$\tuple{x,y}$ is simply $\{\{x\},\{x,y\}\}$, and coding is no longer arithmetisation.

\cite{swierczkowski-finite} quotes \cite{boolos-provability}, who describes his proofs as ``incomplete'' and
``irremediably messy'' (page~16). \'Swierczkowski's proof of the second incompleteness theorem  is certainly
less messy, because he eliminates virtually all arithmetical arguments. The Isabelle/HOL proofs are of
course complete, and represent the first machine-assisted proof of the second incompleteness theorem. The
explicit derivations in the HF calculus are necessarily messy, because they are strings of low-level logical
inferences. But with few exceptions, the statements actually proved are straightforward; generally, they
prove that various coded operations do exactly what they are supposed to do.

The rest of the paper discusses Isabelle/HOL (\ref{sec:background}) and the fundamental definitions
underlying the proofs (\ref{sec:basics}). Techniques used to formalise G\"odel-numbering are briefly
sketched (\ref{sec:coding}). The steps leading to the first incompleteness theorem
is then described (\ref{sec:1st}).
One small but interesting finding concerns the technique for proving the second incompleteness theorem. The
descriptions given by both \cite{boolos-provability} and \cite{swierczkowski-finite} 
are potentially misleading, if not actually wrong (\ref{sec:2nd}). 
Another finding is that {\'S}wierczkowski's proof is actually incomplete, with a
significant gap which I have closed using methods quite different from the ones he outlined (\ref{sec:HF}). 
A brief section concludes the paper (\ref{sec:conclusions}).

Note that this
paper contains no definitions or proofs as conventionally understood in mathematics; rather, it describes
definitions and formal proofs that have been conducted in Isabelle/HOL, and lessons learned from them.
Our focus below concerns such logical issues revealed by the Isabelle/HOL development. Technological aspects
of this development are discussed in a companion paper \cite{paulson-incompl-ar}. 
In order to save space, standard definitions involving the incompleteness theorems are
not presented below except where they need to be discussed specifically. This material is widely available,
and \cite{swierczkowski-finite} can be downloaded from an Internet archive.\footnote{\url{http://
journals.impan.gov.pl/dm/Inf/422-0-1.html}}

\section{Background.} \label{sec:background}
These proofs were conducted using Isabelle/HOL, an interactive theorem prover \citep{isa-tutorial}.
Therefore \emph{all} proofs are conducted in a formal calculus: higher-order logic. Nevertheless, there
is an enormous difference between proofs carried out Isabelle/HOL's native logic and those carried out in a
formal calculus specified within Isabelle/HOL\@. Interactive theorem provers typically hide the underlying
calculus as much as possible through automatic simplifiers and other tools, trying to create the illusion
that the user is writing a rigorous but flexible mathematical document. A logical calculus formalised within
Isabelle/HOL is an inductively defined set, and a proof within this calculus is a demonstration that a
particular object (representing a formula) belongs to that set. Isabelle's automation assists with such
demonstrations, but they are nevertheless long and all but incomprehensible.

Before formalising the logical calculus, we must formalise the syntax of terms and formulas. A crucial
question is the treatment of bound variables. The names of bound variables are typically regarded as
significant, so that $\ex{xy}{x>y}$ and $\ex{vw}{v>w}$ are distinct (albeit logically equivalent) formulas.
With such an approach, renaming a bound variable is an explicit step. G\"odel's proofs make heavy use of
explicit formulas with many quantifiers, and also require induction over the structure of formulas.
Having to rename bound variables complicates proofs considerably.

\emph{Nominal Isabelle} is a formal theory developed within Isabelle/HOL in order to support reasoning about named
bound variables \citep{urban-general}. Variable names are significant where they appear free, but variable
binding constructions are quotiented with respect to the bound variable names, so that $\ex{xy}{x>y}$ and $
\ex{vw}{v>w}$ denote the same formula exactly as $\{0,1\}$ and $\{1,0\}$ denote the same set. Permutations
on names are the key underlying mechanism, for which can be derived the function $\supp(\alpha)$,
which coincides with the set of free variables in $\alpha$ when $\alpha$ is something like a term or formula.
When performing induction on a formula, these mechanisms can ensure that
any bound variables inside the formula are distinct from those of any other formulas that we are interested
in. Thus we can avoid the many problems reported by \cite{oconnor-incompleteness}, who
formalised the first incompleteness theorem using Coq.

One penalty that must be paid in exchange for these advantages is that any function defined on formulas must
use bound variables sensibly (for example, we may not define the set of variables \emph{bound} in a
formula). While the formal definition of ``sensibly'' admits all the definitions required for the
incompleteness theorems, proving this property required specialised skills
(I frequently called upon Christian Urban for assistance), and they can be very demanding of processor time.

For the coding of formulas, bound variables can be formalised using the 
nameless approach of \cite{debruijn72}.
Bound variable occurrences are designated by non-negative integers: 0 for the innermost bound variable
and increasing for each intervening quantifier. Substitution and abstraction can be defined
easily. The main drawback of eliminating bound variable names in this manner is a complete loss of
readability, but that is of no importance for coding. The Isabelle/HOL development proves an
exact correspondence between the syntax of terms and formulas defined using Nominal Isabelle and the codes
of terms and formulas. This correspondence extends to syntactic operations, such as substitution, encoded
using a combination of {\'S}wierczkowski's and de Bruijn's techniques. There is no need to formalise the
nominal theory in the HF calculus, and the complications would be considerable.

A sceptical reader is entitled to ask why we should trust this complicated software and the mysterious
nominal theory. We gain confidence in it---as with all human artefacts---through a combination of personal
experience, its reputation and an understanding of its design. Isabelle/HOL has now been used in a great
many substantial projects by hundreds of users, giving strong reasons to accept that it is a correct
implementation of higher-order logic. The nominal theory is a definitional extension of this logic,
all concepts ultimately reducible to HOL primitives. 
The formally verified correspondence between nominal syntax and de Bruijn syntax, mentioned above, is further
evidence for its correctness. The formal development itself presents a proof of the incompleteness theorems
at a level of detail vastly greater than can be found in any published account. Moreover, this formal
development is a live document: our sceptic can load it into Isabelle/HOL, point to any part of any proof,
and quickly see what has to be proved at that point. Transparency is the best response to scepticism.

\section{The Isabelle/HOL formalisation: fundamentals.} \label{sec:basics}
Let us see what typical definitions and proofs
look like in Isabelle/HOL\@. One claim for
this work is that the machine proofs are readable, at least to a limited extent, allowing this very
lengthy and complicated series of definitions and proofs to be examined.

The hereditarily finite sets are recursively defined as finite sets of hereditarily finite sets.
\cite{swierczkowski-finite} presents a first-order theory having a constant 0 (the empty set), a binary
operation symbol~$\lhd$ (augmentation, or ``eats''),
a relation symbol~$\in$ (membership) as well as equality, satisfying
the following axioms:
\begin{gather}
z=0 \bimp \all{x}{x\not\in z} \tag{HF1}\\
z=x\lhd y \bimp \all{u}{u\in z\bimp u\in x\lor u=y} \tag{HF2}\\
\phi(0) \land \all{xy}{\phi(x)\land\phi(y)\to\phi(x\lhd y)}\to \all{x}{\phi(x)} \tag{HF3}
\end{gather}
The third axiom expresses induction. \cite{swierczkowski-finite} develops the necessary elements of this set
theory, including functions, ordinals (which are simply the natural numbers) and definitional principles.
\cite{kirby-addition} presents an elegant generalisation of ordinal addition to the universe of sets.
Formalising such material in Isabelle/HOL is routine. 

The first milestone in proving the incompleteness theorems is to formalise the syntax of the HF calculus.
Remember, in Isabelle/HOL, mathematics is expressed in higher-order logic. This is a typed formalism, and
the following declaration establishes a recursive type \isa{tm} of HF terms. The type \isa{name} has already
been established, using the nominal framework, as the type of variable names for this calculus.
\begin{isabelle}
\isacommand{nominal\_datatype}\ tm\ =\ Zero\ |\ Var\ name\ |\ Eats\ tm\ tm
\end{isabelle}
This declares that a term is either \isa{Zero} or has the form \isa{Var i}, where \isa{i} is a name, or has
the form \isa{Eats t1 t2} for terms \isa{t1} and \isa{t2}.

It is now possible to define the type \isa{fm} of HF formulas.
\begin{isabelle}
\isacommand{nominal\_datatype}\ fm\ =\isanewline
\ \ \ \ Mem\ tm\ tm\ \ \ \ (\isakeyword{infixr}\ "IN"\ 150)\isanewline
\ \ |\ Eq\ tm\ tm\ \ \ \ \ (\isakeyword{infixr}\ "EQ"\ 150)\isanewline
\ \ |\ Disj\ fm\ fm\ \ \ (\isakeyword{infixr}\ "OR"\ 130)\isanewline
\ \ |\ Neg\ fm\isanewline
\ \ |\ Ex\ x::name\ f::fm\ \isakeyword{binds}\ x\ \isakeyword{in}\ f
\end{isabelle}
The HF calculus includes an existential quantifier, denoted \isa{Ex}, which involves variable binding via
the nominal framework. The \isakeyword{infixr} declarations provide an alternative syntax for the membership
relation, the equality relation, and disjunction. A formula can also be a negation. The other
logical connectives are introduced later as abbreviations.

Substitution is often problematical to formalise, but here it is straightforward. Substitution of a term
\isa{x} for a variable \isa{i} is defined as follows:
\begin{isabelle}
\isacommand{nominal\_primrec}\ subst\ ::\ "name\ \isasymRightarrow \ tm\ \isasymRightarrow \ tm\
\isasymRightarrow \ tm"\isanewline
\ \ \isakeyword{where}\isanewline
\ \ \ "subst\ i\ x\ Zero\ \ \ \ \ \ \ =\ Zero"\isanewline
\ |\ "subst\ i\ x\ (Var\ k)\ \ \ \ =\ (if\ i=k\ then\ x\ else\ Var\ k)"\isanewline
\ |\ "subst\ i\ x\ (Eats\ t\ u)\ =\ Eats\ (subst\ i\ x\ t)\ (subst\ i\ x\ u)"
\end{isabelle}

For substitution within a formula, we normally expect
issues concerning the capture of a bound variable. Note that the result of substituting the term \isa{x} for
the variable \isa{i} in the formula \isa{A} is written \isa{A(i::=x)}.
\begin{isabelle}
\isacommand{nominal\_primrec}\ \ subst\_fm\ ::\ "fm\ \isasymRightarrow \ name\ \isasymRightarrow \ tm\
\isasymRightarrow \ fm"\isanewline
\ \ \isakeyword{where}\isanewline
\ \ \ \ Mem:\ \ "(Mem\ t\ u)(i::=x)\ \ =\ Mem\ (subst\ i\ x\ t)\ (subst\ i\ x\ u)"\isanewline
\ \ |\ Eq:\ \ \ "(Eq\ t\ u)(i::=x)\ \ \ =\ Eq\ \ (subst\ i\ x\ t)\ (subst\ i\ x\ u)"\isanewline
\ \ |\ Disj:\ "(Disj\ A\ B)(i::=x)\ =\ Disj\ (A(i::=x))\ (B(i::=x))"\isanewline
\ \ |\ Neg:\ \ "(Neg\ A)(i::=x)\ \ \ \ =\ Neg\ (A(i::=x))"\isanewline
\ \ |\ Ex:\ \ \ "atom\ j\ \isasymsharp \ (i,\ x)\ \isasymLongrightarrow \ (Ex\ j\ A)(i::=x)\ =\ Ex\ j\
(A(i::=x))"
\end{isabelle}
Substitution is again straightforward in the first four cases (membership, equality, disjunction, negation).
In the existential case, the precondition \isa{atom\ j\ \isasymsharp \ (i,\ x)} (pronounced ``\isa{j} is
fresh for \isa{i} and \isa{x}'') essentially says that \isa{i} and \isa{j} must be different names with
\isa{j} not free in \isa{x}. We do not need to supply a mechanism for renaming the bound variable, as that
is part of the nominal framework, which in most cases will choose a sufficiently fresh bound variable at the
outset. The usual properties of substitution (commutativity, for example) have simple proofs by
induction on formulas. In contrast, \cite{oconnor-phd}  needed to combine three substitution lemmas in a
simultaneous proof by induction, a delicate argument involving 1900 lines of Coq.

The HF proof system is an inductively defined predicate, where \isa{H\ \isasymturnstile \ A} means that the
formula \isa{A} is provable from the set of formulas~\isa{H}\@.
\begin{isabelle}
\isacommand{inductive}\ hfthm\ ::\ "fm\ set\ \isasymRightarrow \ fm\ \isasymRightarrow \ bool"\
(\isakeyword{infixl}\ "\isasymturnstile "\ 55)\isanewline
\ \ \isakeyword{where}\isanewline
\ \ \ \ Hyp:\ \ \ \ "A\ \isasymin \ H\ \isasymLongrightarrow \ H\ \isasymturnstile \ A"\isanewline
\ \ |\ Extra:\ \ "H\ \isasymturnstile \ extra\_axiom"\isanewline
\ \ |\ Bool:\ \ \ "A\ \isasymin \ boolean\_axioms\ \isasymLongrightarrow \ H\ \isasymturnstile \
A"\isanewline
\ \ |\ Eq:\ \ \ \ \ "A\ \isasymin \ equality\_axioms\ \isasymLongrightarrow \ H\ \isasymturnstile \
A"\isanewline
\ \ |\ Spec:\ \ \ "A\ \isasymin \ special\_axioms\ \isasymLongrightarrow \ H\ \isasymturnstile \
A"\isanewline
\ \ |\ HF:\ \ \ \ \ "A\ \isasymin \ HF\_axioms\ \isasymLongrightarrow \ H\ \isasymturnstile \ A"\isanewline
\ \ |\ Ind:\ \ \ \ "A\ \isasymin \ induction\_axioms\ \isasymLongrightarrow \ H\ \isasymturnstile \
A"\isanewline
\ \ |\ MP:\ \ \ \ \ "H\ \isasymturnstile \ A\ IMP\ B\ \isasymLongrightarrow \ H'\ \isasymturnstile \ A\
\isasymLongrightarrow \ H\ \isasymunion \ H'\ \isasymturnstile \ B"\isanewline
\ \ |\ Exists:\ "H\ \isasymturnstile \ A\ IMP\ B\ \isasymLongrightarrow \isanewline
\ \ \ \ \ \ \ \ \ \ \ atom\ i\ \isasymsharp \ B\ \isasymLongrightarrow \ \isasymforall C\isasymin H.\ atom\
i\ \isasymsharp \ C\ \isasymLongrightarrow \ H\ \isasymturnstile \ (Ex\ i\ A)\ IMP\ B"
\end{isabelle}
Note that the existential rule is subject to the condition that the bound variable, \isa{i}, is fresh with
respect to~\isa{B} and the formulas in \isa{H}\@. The definitions of \isa{boolean\_axioms}, etc., are taken from
\cite{swierczkowski-finite}. He formalised a simpler inference system, with theorems of the form \isa{\
\isasymturnstile \ A}. Introducing \isa{H} allows a proof of the deduction theorem and the derivation of
a sort of sequent calculus, a practical necessity if we are to conduct proofs in this formal calculus.

Another deviation from \cite{swierczkowski-finite} is the inclusion of \isa{extra\_axiom}. It is a
parameter of the entire development; it can be any formula that is true under
the Tarski truth-definition.\footnote{This is formalised as the function \isa{eval\_fm}, which is 
presented in the companion paper \citep[section 3.1]{paulson-incompl-ar}.
The constraint that \isa{extra\_axiom} must be true is not shown here.}
Its purpose is to generalise the statements of the incompleteness
theorems, which {\'S}wierczkowski proved only for one specific calculus. \cite{oconnor-incompleteness} has
gone further to prove the first incompleteness theorem even for infinite extensions of the calculus.

The incompleteness theorems require the definition of a great many predicates, mostly for coding the syntax
of terms and formulas, and operations on them. It may be instructive to look at a very simple definition,
namely of the subset relation:
\begin{isabelle}
\isacommand{nominal\_primrec}\ Subset\ ::\ "tm\ \isasymRightarrow \ tm\ \isasymRightarrow \ fm"\ \ \
(\isakeyword{infixr}\ "SUBS"\ 150)\isanewline
\ \ \isakeyword{where}\ "atom\ z\ \isasymsharp \ (t,\ u)\ \isasymLongrightarrow \ t\ SUBS\ u\ =\ All2\ z\ t\
((Var\ z)\ IN\ u)"
\end{isabelle}
This introduces \isa{SUBS} as the name of the subset relation, which is defined using a bounded quantifier
by $t\subseteq u \iff \all{(z\in t)}{z\in u}$. Note that \isa{All2} is our syntax for a bounded universal
quantifier. The condition \isa{atom\ z\ \isasymsharp \ (t,u)} states that the quantified variable (\isa{z})
must be fresh for the terms \isa{t} and \isa{u}. In other words, and in contrast to some treatments, the
bound variable is a parameter of the definition rather than being fixed; however, the choice of
\isa{z} cannot affect the denotation of the right-hand side, thanks to quotienting.

Proving the elementary properties of the subset relation within the HF calculus is extremely tedious, over
200 lines of proof script. Extensionality must be proved by induction within the calculus:
\begin{isabelle}
\isacommand{lemma}\ Extensionality:\ "H\ \isasymturnstile \ x\ EQ\ y\ IFF\ (x\ SUBS\ y\ AND\ y\ SUBS\ x)"
\end{isabelle}
The length of these trivial proofs might be taken as a
sign that mechanising the incompleteness theorems is infeasible. It is fortunate that proofs of apparently
more advanced properties do not get longer and longer, even when we come to prove the Hilbert-Bernays
derivability conditions.

\cite{swierczkowski-finite} discusses
$\Sigma$ formulas, constructed from atomic formulas using conjunction, disjunction, existential
quantification and bounded universal quantification. \emph{Strict} $\Sigma$ formulas contain no terms
other than variables, and the bound~$j$ in $\forall{(i\in j)}\,A$ must not be free in the quantified body, $A$.
\begin{isabelle}
\isacommand{inductive}\ ss\_fm\ ::\ "fm\ \isasymRightarrow \ bool"\ \isakeyword{where}\isanewline
\ \ \ \ MemI:\ \ "ss\_fm\ (Var\ i\ IN\ Var\ j)"\isanewline
\ \ |\ DisjI:\ "ss\_fm\ A\ \isasymLongrightarrow \ ss\_fm\ B\ \isasymLongrightarrow \ ss\_fm\ (A\ OR\
B)"\isanewline
\ \ |\ ConjI:\ "ss\_fm\ A\ \isasymLongrightarrow \ ss\_fm\ B\ \isasymLongrightarrow \ ss\_fm\ (A\ AND\
B)"\isanewline
\ \ |\ ExI:\ \ \ "ss\_fm\ A\ \isasymLongrightarrow \ ss\_fm\ (Ex\ i\ A)"\isanewline
\ \ |\ All2I:\ "ss\_fm\ A\ \isasymLongrightarrow \ atom\ j\ \isasymsharp \ (i,A)\ \isasymLongrightarrow \ ss\_fm\ (All2\ i\ (Var\ j)\ A)"
\end{isabelle}
One advantage of formal proof is that these conditions are immediately evident, when they may not be clear
from an informal presentation. \cite{swierczkowski-finite} does not impose the last condition (on the bound
of a universal quantifier), but it greatly simplifies the main induction needed to reach the second
incompleteness theorem. (If we are only interested in formalising the first incompleteness theorem, we can
use a more generous notion of $\Sigma$ formula, allowing atomic formulas and their negations over arbitrary
terms.) Formally, a $\Sigma$ formula is defined to be any formula that can be proved equivalent (in the HF
calculus) to a strict $\Sigma$ formula:
\begin{isabelle}
\ \ "Sigma\_fm\ A\ \isasymlongleftrightarrow \ (\isasymexists B.\ ss\_fm\ B\ \&\ supp\ B\ \isasymsubseteq \
supp\ A\ \&\ \isacharbraceleft \isacharbraceright \ \isasymturnstile \ A\ IFF\ B)"
\end{isabelle}
The condition \isa{supp\ B\ \isasymsubseteq \ supp\ A} essentially means that every variable free in \isa{B}
must also be free in \isa{A}\@. After a certain amount of effort, it is possible to derive the expected
properties of $\Sigma$ formulas and ultimately to reach a key result based on this concept:
\begin{isabelle}
\isacommand{theorem}\ Sigma\_fm\_imp\_thm:\ "\isasymlbrakk Sigma\_fm\ A;\ ground\_fm\ A;\ eval\_fm\ e0\ A\isasymrbrakk
\ \isasymLongrightarrow \ \isacharbraceleft \isacharbraceright \ \isasymturnstile \ A"
\end{isabelle}
If \isa{A} is a true $\Sigma$ sentence, then \isa{\isasymturnstile A}. This result reduces the task of proving
\isa{\isasymturnstile \ A} in the formal calculus to proving that \isa{A} holds (written \isa{eval\_fm\ e0\
A}) in Isabelle/HOL's native higher-order logic.

\section{The Isabelle/HOL formalisation: The coding of syntax.} \label{sec:coding}
The coding of terms,
formulas, substitution, the HF axioms and ultimately the provability predicate is straightforward to formalise.
\cite{goedel31} and \cite{swierczkowski-finite} present full details. Many other authors prefer
to simplify matters via repeated appeals to Church's thesis. Even the detailed presentations
mentioned above omit any demonstration that the definitions are correct. The
\emph{proof formalisation condition} for the
provability predicate (written \isa{PfP} below) is typically stated with a minimum of justification:
\begin{isabelle}
\isacommand{theorem}\ proved\_iff\_proved\_Pf:\ "\isacharbraceleft \isacharbraceright
\ \isasymturnstile \ A\ \isasymlongleftrightarrow \ \isacharbraceleft \isacharbraceright
\ \isasymturnstile \ PfP\ \isasymlceil A\isasymrceil "
\end{isabelle}
One could argue that there is no need for the definitions to be correct in every detail, provided they convince
the reader that correct and suitable definitions exist.
However, only correct definitions can be verified in Isabelle/HOL\@. 
Most of these proofs are indeed routine,
though in places (for example, in the specification of an instance of the HF induction axiom) extremely tedious.

The \cite{debruijn72} representation of variable binding requires new versions of the syntactic predicates
for ``formula'', ``substitution'', etc. The coding of terms and formulas is done by first translating them
from nominal syntax to de Bruijn syntax. In verifying the coding predicates, we also verify this
translation.

A standard treatment of de Bruijn syntax requires defining two operations: abstraction and substitution.
Abstraction replaces free occurrences of a given term by a new bound variable, represented by a numeric
index; the resulting formula is ill-formed until a matching quantifier is prefixed to it. Substitution is
the inverse of abstraction, replacing the outermost bound variable (after a quantifier has been stripped from
a formula) by some given term. For the incompleteness theorems, both operations can be simplified:
abstraction replaces a free variable by a bound variable, and substitution replaces a free variable by a
given term. Abstraction is needed to formalise the construction of a formula, because it is a
necessary step before a quantifier can be attached.

The interplay of these various points can be seen below:
\begin{isabelle}
\isacommand{definition}\ MakeForm\ ::\ "hf\ \isasymRightarrow \ hf\ \isasymRightarrow \ hf\ \isasymRightarrow \ bool"\isanewline
\isakeyword{where}\ "MakeForm\ y\ u\ w\ \isasymequiv \isanewline
\ \ \ \ \ y\ =\ q\_Disj\ u\ w\ \isasymor \ y\ =\ q\_Neg\ u\ \isasymor \isanewline
\ \ \ \ \ (\isasymexists v\ u'.\ AbstForm\ v\ 0\ u\ u'\ \isasymand \ y\ =\ q\_Ex\ u')"
\end{isabelle}
Thus \isa{y} is the code of a formula constructed from existing formulas \isa{u} and \isa{v}
provided \isa{y} codes the disjunction $\isa{u} \lor \isa{v}$, the negation $\neg\isa{u}$ or the
existential formula $\exists(\isa{u'})$, where \isa{u'} has been obtained by abstracting \isa{u} over some variable, \isa{v}.
The predicate \isa{AbstForm} performs de Bruijn abstraction over a formula;
its definition is complicated, and omitted here.
Note that the codes of quantified formulas do not mention the names of bound variables.

This predicate is given by a higher-order logic formula, and therefore at the level of the meta-theory.
Working at this level eliminates the need to construct HF proofs, and most of the correctness properties we
need can be proved in this manner. However, in order to perform the diagonalisation argument and exhibit the
undecidable formula, we need a version of every coding predicate as an HF formula. Therefore, each predicate
must be defined on both levels:
\begin{isabelle}
\isacommand{nominal\_primrec}\ MakeFormP\ ::\ "tm\ \isasymRightarrow \ tm\ \isasymRightarrow \ tm\ \isasymRightarrow \ fm"\isanewline
\isakeyword{where}\ "\isasymlbrakk atom\ v\ \isasymsharp \ (y,u,w,au);\ atom\ au\ \isasymsharp \ (y,u,w)\isasymrbrakk \ \isasymLongrightarrow \isanewline
\ \ MakeFormP\ y\ u\ w\ =\isanewline
\ \ \ \ y\ EQ\ Q\_Disj\ u\ w\ OR\ y\ EQ\ Q\_Neg\ u\ OR\isanewline
\ \ \ \ Ex\ v\ (Ex\ au\ (AbstFormP\ (Var\ v)\ Zero\ u\ (Var\ au)\ AND\ y\ EQ\ Q\_Ex\ (Var\ au)))"
\end{isabelle}
As we saw above in the definition of \isa{Subset}, constraints are required on all quantified variables.
Here there are only two, but to define \isa{AbstForm} requires 12 bound variables.
The necessary declarations are lengthy and messy, and put a heavy burden on the nominal package
(proofs run very slowly), but the alternative of having to rename explicit bound variables is also unattractive.

\section{The Isabelle/HOL formalisation: first incompleteness theorem.} \label{sec:1st}
The diagonalisation theorem
is now easily reached. Continuing to follow \cite{swierczkowski-finite}, the next step is to define
a function $K$ such that $\vdash K(\quot{\phi}) =\quot{\phi(\quot{\phi})}$. Formally, $K$
is a \emph{pseudo-function}, represented by the single-valued relation \isa{KRP}, taking two arguments.
The following result is not difficult to obtain,
given the existing coding of substitution, and some other steps that will be discussed later.
This theorem does not require a proof within the HF calculus, but follows from \isa{Sigma\_fm\_imp\_thm}
because it is a sentence (coded syntax contains no free variables) and a $\Sigma$ formula.
\begin{isabelle}
\isacommand{lemma}\ prove\_KRP:\ "\isacharbraceleft \isacharbraceright
\ \isasymturnstile \ KRP\ \isasymlceil Var\ i\isasymrceil \ \isasymlceil A\isasymrceil \ \isasymlceil A(i::=\isasymlceil A\isasymrceil )\isasymrceil "
\end{isabelle}

The property of being single-valued is easily stated, but it is neither a sentence nor a $\Sigma$ formula.
Proving this result requires about 600 lines of explicit reasoning steps in the HF calculus, verifying
that substitution over terms or formulas yields a unique result.
\begin{isabelle}
\isacommand{lemma}\ KRP\_unique:\ "\isacharbraceleft KRP\ v\ x\ y,\ KRP\ v\ x\ y'\isacharbraceright \ \isasymturnstile \ y'\ EQ\ y"
\end{isabelle}
The diagonal lemma is now reached by the standard argument. 
It concerns an arbitrary formula, \isa{\isasymalpha}, presumably containing \isa{i} as a free variable.
Note that \isa{\isasymalpha (i::=\isasymlceil \isasymdelta \isasymrceil)} denotes the result of replacing 
\isa{i} by \isa{\isasymlceil \isasymdelta \isasymrceil}.
The \isakeyword{obtains} syntax represents
a form of existential quantification, here asserting the existence of an HF formula \isa{\isasymdelta}
satisfying the two properties shown.
\begin{isabelle}
\isacommand{lemma}\ diagonal:\ \isanewline
\ \ \isakeyword{obtains}\ \isasymdelta \ \isakeyword{where}\ "\isacharbraceleft \isacharbraceright
\ \isasymturnstile \ \isasymdelta \ IFF\ \isasymalpha (i::=\isasymlceil \isasymdelta \isasymrceil )"\ \ "supp\ \isasymdelta \ =\ supp\ \isasymalpha \ -\ \isacharbraceleft atom\ i\isacharbraceright "
\end{isabelle}
The second part of the conclusion, namely \isa{supp\ \isasymdelta \ =\ supp\ \isasymalpha \ -\ \isacharbraceleft atom\ i\isacharbraceright},
states that the free variables of the formula \isa{\isasymdelta} are those of \isa{\isasymalpha}
with the exception of~\isa{i}; it is necessary in order to show that the undecidable formula
is actually a sentence.

The first incompleteness theorem itself can now be proved. Figure~\ref{fig:Goedel-I} presents the full text.
Even a reader who is wholly unfamiliar with Isabelle/HOL should be able to see something intelligible
in this proof script. Assuming consistency of the calculus, formalised as
\hbox{\isa{\isasymnot \ \isacharbraceleft \isacharbraceright \ \isasymturnstile \ Fls}}
(falsity is not provable), we obtain a formula \isa{\isasymdelta} satisfying the properties shown,
in particular \isa{\isasymnot \ \isacharbraceleft \isacharbraceright \ \isasymturnstile \ \isasymdelta}
and  \isa{\isasymnot \ \isacharbraceleft \isacharbraceright \ \isasymturnstile \ Neg\ \isasymdelta}.
Lines beginning with commands such as \isacommand{obtain}, \isacommand{hence},
\isacommand{show} introduce assertions to be proved. The details of the reasoning
may be unclear, but milestones such as
\isa{"\isacharbraceleft \isacharbraceright \ \isasymturnstile \ \isasymdelta \ IFF\ Neg\ (PfP\ \isasymlceil \isasymdelta \isasymrceil )"} and
\hbox{\isa{"\isasymnot \ \isacharbraceleft \isacharbraceright \ \isasymturnstile \ \isasymdelta "}}
are visible, as references to previous named results. This legibility, however limited, is possible because
the entire Isabelle/HOL proof is written in the structured Isar language
\citep{wenzel-isabelle/isar}. Only the HF calculus proofs remain unintelligible:
it is not easy to impose structure on those.

\begin{figure}
\begin{center}
\begin{isabelle}
\isacommand{theorem}\ Goedel\_I:\isanewline
\ \ \isakeyword{assumes}\ "\isasymnot \ \isacharbraceleft \isacharbraceright \ \isasymturnstile \ Fls"\isanewline
\ \ \isakeyword{obtains}\ \isasymdelta \ \isakeyword{where}\ "\isacharbraceleft \isacharbraceright
\ \isasymturnstile \ \isasymdelta \ IFF\ Neg\ (PfP\ \isasymlceil \isasymdelta \isasymrceil )"
\ \ "\isasymnot \ \isacharbraceleft \isacharbraceright \ \isasymturnstile \ \isasymdelta "\ \ "\isasymnot \ \isacharbraceleft \isacharbraceright \ \isasymturnstile \ Neg\ \isasymdelta "\isanewline
\ \ \ \ \ \ \ \ \ \ \ \ \ \ \ \ \ "eval\_fm\ e\ \isasymdelta "\ \ "ground\_fm\ \isasymdelta "\isanewline
\isacommand{proof}\ -\isanewline
\ \ \isacommand{obtain}\ \isasymdelta \ \isakeyword{where}\ \ \ \ \ \ \ \ \ "\isacharbraceleft \isacharbraceright \ \isasymturnstile \ \isasymdelta \ IFF\ Neg\ ((PfP\ (Var\ i))(i::=\isasymlceil \isasymdelta \isasymrceil ))"\isanewline
\ \ \ \ \ \ \ \ \ \ \ \ \isakeyword{and}\ [simp]:\ "supp\ \isasymdelta \ =\ supp\ (Neg\ (PfP\ (Var\ i)))\ -\ \isacharbraceleft atom\ i\isacharbraceright "\isanewline
\ \ \ \ \isacommand{by}\ (metis\ SyntaxN.Neg\ diagonal)\isanewline
\ \ \isacommand{hence}\ diag:\ "\isacharbraceleft \isacharbraceright \ \isasymturnstile \ \isasymdelta \ IFF\ Neg\ (PfP\ \isasymlceil \isasymdelta \isasymrceil )"\isanewline
\ \ \ \ \isacommand{by}\ simp\isanewline
\ \ \isacommand{hence}\ np:\ "\isasymnot \ \isacharbraceleft \isacharbraceright \ \isasymturnstile \ \isasymdelta "\isanewline
\ \ \ \ \isacommand{by}\ (metis\ Con\ Iff\_MP\_same\ Neg\_D\ proved\_iff\_proved\_Pf)\isanewline
\ \ \isacommand{hence}\ npn:\ "\isasymnot \ \isacharbraceleft \isacharbraceright \ \isasymturnstile \ Neg\ \isasymdelta "\ \isacommand{using}\ diag\isanewline
\ \ \ \ \isacommand{by}\ (metis\ Iff\_MP\_same\ NegNeg\_D\ Neg\_cong\ proved\_iff\_proved\_Pf)\isanewline
\ \ \isacommand{moreover}\ \isacommand{have}\ "eval\_fm\ e\ \isasymdelta "\ \isacommand{using}\ hfthm\_sound\ [\isakeyword{where}\ e=e,\ OF\ diag]\isanewline
\ \ \ \ \isacommand{by}\ simp\ (metis\ Pf\_quot\_imp\_is\_proved\ np)\isanewline
\ \ \isacommand{moreover}\ \isacommand{have}\ "ground\_fm\ \isasymdelta "\ \isanewline
\ \ \ \ \isacommand{by}\ (auto\ simp:\ ground\_fm\_aux\_def)\isanewline
\ \ \isacommand{ultimately}\ \isacommand{show}\ ?thesis\isanewline
\ \ \ \ \isacommand{by}\ (metis\ diag\ np\ npn\ that)\isanewline
\isacommand{qed}
\end{isabelle}

\caption{Proof of the first incompleteness theorem}
\label{fig:Goedel-I}
\end{center}
\end{figure}

\section{Issues involving the second incompleteness theorem.} \label{sec:2nd}
My object in writing this paper is not
to discuss the formalisation in general, but to examine the specific consequences
of basing the development on HF set theory rather than Peano arithmetic. A further aim is to
look at a crucial step in the proof of the second incompleteness theorem that is not described
especially well in other presentations.

It is well-known that the theorem follows easily from the Hilbert-Bernays
derivability conditions \cite[p.\ts15]{boolos-provability}, one of which is $\vdash \Pf(\quot{\phi}) \to \Pf(\quot{\Pf(\quot{\phi})})$.
This result is a consequence of the theorem
\begin{gather}
\text{if $\alpha$ is a $\Sigma$ sentence, then } \vdash \alpha\to \Pf(\quot{\alpha}), \label{eqn:star}
\end{gather}
which can be proved by a tricky induction on the construction of $\alpha$ as a strict $\Sigma$ formula.

For this proof, the system of coding is extended to allow variables in codes. If we regard variables
as indexed by positive integers, then the variable $x_i$ is normally coded by the term $\SUCC^i (0)$,
where $\SUCC(x) = x\lhd x$ is the usual successor function. Similarly, the formula
$x_1=x_2$ is normally coded by the term $\tuple{\quot{=},\quot{x_1},\quot{x_2}}$.
If variables are preserved rather than coded, we instead get the term  $\langle\quot{=},x_1,x_2\rangle$.
In general, $\vquot{\alpha}_V$ designates the coding of $\alpha$ where all variables from the set~$V$
are preserved as variables in the result, while all other variables are coded by constant terms.
\cite{swierczkowski-finite} calls this \emph{pseudo-coding}.

Imagine that we could define in HF a function $Q$ such that
\begin{align}
Q(0) & = \quot{0} = 0
	\label{eqn:Q0} \\
Q(x\lhd y) &= \tuple{\quot{\lhd},Q(x),Q(y)}
	\label{eqn:Q1}
\end{align}
Then we would have $Q(x) = \quot{t}$, where $t$ is some canonical term denoting the set~$x$.
[\cite{swierczkowski-finite} introduces a total ordering on HF to make this possible, as
discussed below.]
Suppose that $\alpha$ is a formula whose set of free variables is $V=\{x_1,\ldots,x_n\}$.
Given the theorem $\vdash \alpha$, obtain $\vdash \Pf(\quot{\alpha})$ by the proof formalisation condition,
then successively replace $x_i$ by $Q(x_i)$, for $i=1$, \ldots,~$n$. The replacements are possible
because the HF calculus includes a rule to substitute a term $t$ for a variable~$x$ in the formula~$\phi$:
\[ \frac{H\vdash \phi}{H\vdash \phi(x/t)} \]
Performing the replacements requires the analogue of this substitution rule as encoded in
the provability predicate, $\Pf$. For example, we can obtain the following series of theorems:
\begin{alignat*}{2}
& \vdash y\in (x \lhd y) \\
& \vdash \Pf \quot{y\in (x \lhd y)}  &\qquad& \text{proof formalisation condition}\\
& \vdash \Pf \tuple{\quot{\in},\quot{y},\tuple{\quot{\lhd},\quot{x},\quot{y}}}  && \text{definition of coding}\\
& \vdash \Pf \tuple{\quot{\in},\quot{y},\tuple{\quot{\lhd},Q(x),\quot{y}}}  && \text{replacement of $x$}\\
& \vdash \Pf \tuple{\quot{\in},Q(y),\tuple{\quot{\lhd},Q(x),Q(y)}}  && \text{replacement of $y$}
\end{alignat*}
To simplify the notation, let $\alpha (V/Q)$ abbreviate $\alpha(x_1/Q(x_1),\ldots,x_n/Q(x_n))$,
the result of simultaneously replacing every free variable $x_i$ in  $\alpha$ by the term~$Q(x_i)$.
As a further simplification, let us write $\vquot{t}_V \,(Q)$  instead of $\vquot{t}_V \,(V/Q)$.
Then the sequence of steps above can also be written
\begin{alignat*}{2}
& \vdash \Pf \quot{y\in (x \lhd y)}  &\qquad&\\
& \vdash \Pf \vquot{y\in (x \lhd y)}_{\{x\}} \,(Q) && \text{replacement of $x$}\\
& \vdash \Pf \vquot{y\in (x \lhd y)}_{\{x,y\}} \,(Q)  && \text{replacement of $y$}
\end{alignat*}

A crucial part of the reasoning is that the replacement of $\quot{y}$ by $Q(y)$ leaves
the occurrences of $Q(x)$ unchanged. That holds because $Q(x)$ is always the code of a \emph{constant} term,
as can trivially be proved from (\ref{eqn:Q0}) and (\ref{eqn:Q1}) by induction on~$x$. A constant term
is unaffected by substitutions.

{\sloppy
The difficulty with this sketch is that no function $Q(x)$ can exist, because the HF language
has only one function symbol,~$\lhd$. Extending this language with the function symbol $Q$ would require
redoing all the coding and syntactic functions; $Q$ would also need to encode references to itself.
Instead, $Q(x)$ is typically regarded as a ``pseudo-function'': 
it must be defined in the form of a relation $\QR(x,y)$ for which
$\all{x}{\exists! y\, {\QR(x,y)}}$ can be proved.  We must modify the transformations above accordingly.
\cite{boolos-provability} and \cite{swierczkowski-finite} both state that the formula $A(Q(x))$
is an abbreviation for $\ex{x'}{\QR(x,x') \land A(x')}$; the latter author describes a detailed procedure
for replacing occurrences of pseudo-functions from the inside out \cite[p.\ts47]{swierczkowski-finite}.
This suggests the following modified sequence:
\begin{alignat*}{2}
& \vdash \Pf \quot{y\in (x \lhd y)}  &\qquad& \\
& \vdash \Pf \vquot{\ex{x'}{\QR(x,x') \land y\in (x' \lhd y)}}_{\{x\}} && \text{replacement of $x$}\\
& \vdash \Pf \vquot{\ex{y'}{\QR(y,y') \land \ex{x'}{\QR(x,x') \land y'\in (x' \lhd y')}}}_{\{x,y\}}  && \text{replacement of $y$}
\end{alignat*}
Further evidence that this is the intended transformation is the remark \cite[p.\ts45]{boolos-provability}
that the transformed formula, $\Pf(\vquot{\alpha}_V\,(Q))$ in our notation,
``has the same variables free as'' the original formula,~$\alpha $. The difficulty is that this modified
sequence does not work, and neither can any other that leaves the original variables free
in the transformed formula. The explanation is simple: these variables (in particular $x$ above)
range over all values, including the codes of all possible formulas. There is no reason why
$\QR(x,x')$ should be left unchanged after the substitution for~$y$: there is nothing to exclude the
possibility that $x=\quot{y}$, for example. One could argue that the remarks and explanations
that I have cited are true in spirit if not in fact, but they are misleading. We even see
a detailed proof that $Q(x_i)$ is correctly substituted for $x_i$ with reference to the definitions of
the syntactic substitution predicates \cite[p.\ts25]{swierczkowski-finite}, but there is no such term as $Q(x)$.
}

The correct sequence of steps introduces new free variables in the transformed formula,
while simultaneously constraining them as constants on the left-hand side of the $\vdash$ symbol.
\begin{alignat*}{2}
    & \vdash \Pf \quot{y\in x \lhd y}  &\qquad& \\
 \QR(x,x') & \vdash \Pf \vquot{y\in x' \lhd y}_{\{x'\}} && \text{replacement of $x$}\\
 \QR(y,y'), \QR(x,x') & \vdash \Pf \vquot{y'\in x' \lhd y'}_{\{x',y'\}} && \text{replacement of $y$}
\end{alignat*}
Here, $x$ is replaced by~$x'$, constrained by the new assumption $\QR(x,x')$ and then
$y$ is replaced by~$y'$. Now $x'$ is unaffected by substitutions because (given the assumption $\QR(x,x')$)
it can be shown to contain no variables. This reasoning is straightforward enough to conduct formally
in the HF calculus.

This may seem to be a small detail, but as can be seen, it is not difficult to explain correctly.
One could argue that the correct version is actually simpler to explain than the traditional version
involving the pseudo-function $Q$: the notation $\vquot{\alpha}_V\,(Q)$ is no longer necessary.
Eliminating the pseudo-functions from the presentation actually simplifies it.

\section{Issues connected with the use of HF sets.} \label{sec:HF}
The motivation for using hereditarily finite sets rather than Peano arithmetic is that it allows
more natural and simpler proofs. But it appears to complicate
the definition of the function $Q(x)$ mentioned above, which is needed to prove both 
incompleteness theorems. In PA, the analogous function $Z(n)$ is trivial to define \citep[p.\ts165]{goedel-I}: 
there is only one way to write a natural number in the form $\SUCC^n (0)$.

\cite{swierczkowski-finite} eliminates the ambiguity implicit in (\ref{eqn:Q1}) above
by appealing to a total ordering, $<$, on the HF universe. The difficulty is how to define this ordering
within the HF calculus. {\'S}wierczkowski develops the theory, including a definition by recursion 
on the rank of a set, but it does not look easy to formalise in HF\@. Another approach is to define
the function $f:\text{HF}\to \mathbb{N}$ such that $f(x)=\sum\,\{2^{f(y)}\mid y\in x\}$. 
Then we can define $x<y\iff f(x)<f(y)$. Again, the effort to formalise
this theory in HF may be simpler than that needed to formalise the Chinese remainder theorem,
but it is still considerable.

The alternative is to eliminate the need for this ordering. {\'S}wierczkowski has already
completed part of this task. In his proof of the first incompleteness theorem, he introduces
a function $H$ such that $H(\quot{\phi}) = \quot{\quot{\phi}}$. This function is recursively defined
on valid codes, that is, on terms recursively built over natural numbers using ordered pairing.
In fact, $H$ is identical to $Q$ but with a restricted domain,  ensuring that it can 
easily be proved to be a function.

For the second incompleteness theorem, the solution to our conundrum is again to focus on
the corresponding relation, $\QR$. There is no need to prove that this relation describes a function.
All that is necessary in order to prove~(\ref{eqn:star}) is the property
\begin{gather}
\QR(x,x'), \QR(y,y') \vdash x\in y \to \Pf \vquot{x'\in y'}_{\{x',y'\}}, \label{eqn:star-in}
\end{gather}
{\'S}wierczkowski shows that this follows from the lemma
\begin{gather}
\QR(x,x'), \QR(y,y') \vdash x=y \to \Pf \vquot{x'=y'}_{\{x',y'\}}, \label{eqn:star-eq}
\end{gather}
which clearly holds even if $\QR$ does not describe a functional relationship.
A way to prove both (\ref{eqn:star-in}) and (\ref{eqn:star-eq})
can be seen from the following elementary set-theoretic equivalences,
which connect the relations $\in$, $\subseteq$ and~$=$:
\begin{align*}
z\in\emptyset &\iff \bot\\
z\in x\lhd y  &\iff z\in x \lor z=y \\
\emptyset\subseteq z &\iff \top\\
x\lhd y \subseteq z &\iff x \subseteq z \land y\in z \\
x=y &\iff x \subseteq y \land y\subseteq x
\end{align*}

The point of all this is that (\ref{eqn:star-in}) and (\ref{eqn:star-eq}) can be proved by a simultaneous induction:
\begin{gather*}
\QR(x,x'), \QR(y,y') \,\vdash\, (x\in y \to \Pf \vquot{x'\in y'}_{\{x',y'\}}) \land (x\subseteq y \to \Pf \vquot{x'\subseteq y'}_{\{x',y'\}})
\end{gather*}
The induction is on the sum of the lengths of the derivations of $\QR(x,x')$ and $\QR(y,y')$.
Like most of the syntactic predicates used in the incompleteness theorems, $\QR(x,x')$ is defined to hold
provided there exist $k$ and $s$ such that $s$ is a $k$-element sequence representing the conditions
(\ref{eqn:Q0}) and (\ref{eqn:Q1}). Induction on the sum of the lengths allows us to prove
\[ x\in y \to \Pf \vquot{x'\in y'}_{\{x',y'\}} \]
by case analysis on the form of $y$, while proving
\[ x \subseteq y \to \Pf \vquot{x' \subseteq y'}_{\{x',y'\}} \]
by case analysis on the form of $x$.  One case of the reasoning is as follows:
\begin{align*}
x_1\lhd x_2 \subseteq y & \iff x_1 \subseteq y \,\land\, x_2 \in y \\
   & \;\,\Longrightarrow\; \Pf \vquot{x_1'\subseteq y'}_{\{x_1',y'\}} \,\land\, \Pf \vquot{x_2'\in y'}_{\{x_2',y'\}} \\
   & \iff \Pf \vquot{x_1'\lhd x_2'\subseteq y'}_{\{x_1',x_2',y'\}} 
\end{align*}
The formalisation of the entire mutually inductive argument
in the HF calculus requires under 450 lines of Isabelle/HOL\@. The need to define an ordering
on the HF universe has disappeared.

The mechanised proof requires only the simplest induction principles throughout. The basic principle
of the hereditarily finite sets (HF3) is used eight times, mostly to develop the fundamentals
of the HF set theory itself. Complete induction on the natural numbers is used ten times,
while ordinary mathematical induction is used eleven times.
No other form of induction is necessary. \cite{swierczkowski-finite} frequently sketches proofs by induction
on terms or formulas. He suggests induction on the HF ordering, $<$, to prove (\ref{eqn:star-in}) above
and also to prove the bounded quantifier case of the main theorem:
\[\vdash \forall{(j\in i)}\,{\alpha(j)} \to \Pf(\vquot{\forall{(j'\in i)}\,{\alpha(j')}}) \]
Each of these theorems concerns syntactic predicates defined by the existence of a
$k$-element sequence, and is more directly proved by complete induction on~$k$, or rarely
(where there are two sequences, as above) on the sum~$k_1+k_2$.

\section{Discussion and conclusions.} \label{sec:conclusions}
The first mechanised formalisation of G\"odel's (first) incompleteness theorem is due to \cite{shankar-phd}.
It was an astonishing accomplishment given the technology of the 1980s. An interesting technical note
is that \cite{shankar-deBruijn} found de Bruijn indices indispensable in a companion proof (of the Church-Rosser theorem),
but not in his formalisation of the logical calculus. He also used HF set theory, but using
a different axiom system \cite[p.\ts12]{shankar94} that he attributes to Cohen. 
Nineteen years later, \cite{oconnor-incompleteness} mechanised the first theorem using quite different 
methods and the Coq proof assistant. Another proof, by John Harrison,
can be downloaded with his HOL Light proof assistant, \url{http://code.google.com/p/hol-light/}.
There appears to exist no other machine proof of the second incompleteness theorem.

The mechanised incompleteness theorems described above were difficult chiefly because of their sheer size,
and because of the presentational issues discussed from~\ref{sec:2nd} onwards, which resulted in
a great deal of wasted work. But we now have a complete, transparent and machine-checked formalisation
of these landmark results.

\paragraph*{Acknowledgement}
Jesse Alama drew my attention to \cite{swierczkowski-finite}, which was the source material for this project.
Christian Urban assisted with some proofs and wrote some code involving his nominal package. Brian Huffman
assisted with the formalisation of the HF sets. Dana Scott offered advice and drew my attention useful
related work, for example \cite{kirby-addition}.  Matt Kaufmann made insightful comments
on a draft of this paper. The referee made a great many constructive remarks.

\bibliographystyle{rsl}
\bibliography{string,atp,funprog,general,isabelle,theory,crossref}

\begin{thebibliography}{}

\bibitem[\protect\citeauthoryear{Boolos}{Boolos}{1993}]{boolos-provability}
Boolos, G.~S. (1993).
\newblock {\em The Logic of Provability}.
\newblock Cambridge University Press.

\bibitem[\protect\citeauthoryear{de~Bruijn}{de~Bruijn}{1972}]{debruijn72}
de~Bruijn, N.~G. (1972).
\newblock Lambda calculus notation with nameless dummies, a tool for automatic
  formula manipulation, with application to the {Church-Rosser Theorem}.
\newblock {\em Indagationes Mathematicae\/}~{\bf 34}, 381--392.

\bibitem[\protect\citeauthoryear{Feferman}{Feferman}{1986}]{goedel-I}
Feferman, S., editor (1986).
\newblock {\em {Kurt G\"odel}: Collected Works}, Volume~I.
\newblock Oxford University Press.

\bibitem[\protect\citeauthoryear{Franz{\'e}n}{Franz{\'e}n}{2005}]{franzen-guide}
Franz{\'e}n, T. (2005).
\newblock {\em {G\"odel's} Theorem: An Incomplete Guide to Its Use and Abuse}.
\newblock A K Peters.

\bibitem[\protect\citeauthoryear{G{\"o}del}{G{\"o}del}{1931}]{goedel31}
G{\"o}del, K. (1931).
\newblock {\"Uber} formal unentscheidbare {S{\"a}tze} der {Principia
  Mathematica} und verwandter {Systeme I}.
\newblock {\em Monatshefte f{\"u}r Mathematik und Physik\/}~{\bf 38\/}(1),
  173--198.

\bibitem[\protect\citeauthoryear{Kirby}{Kirby}{2007}]{kirby-addition}
Kirby, L. (2007).
\newblock Addition and multiplication of sets.
\newblock {\em Mathematical Logic Quarterly\/}~{\bf 53\/}(1), 52--65.

\bibitem[\protect\citeauthoryear{Nipkow, Paulson, \& Wenzel}{Nipkow
  et~al.}{2002}]{isa-tutorial}
Nipkow, T., Paulson, L.~C., \& Wenzel, M. (2002).
\newblock {\em Isabelle/HOL: A Proof Assistant for Higher-Order Logic}.
\newblock Springer.
\newblock Online at
  \url{http://isabelle.in.tum.de/dist/Isabelle/doc/tutorial.pdf}.

\bibitem[\protect\citeauthoryear{O'Connor}{O'Connor}{2005}]{oconnor-incompleteness}
O'Connor, R. (2005).
\newblock Essential incompleteness of arithmetic verified by {Coq}.
\newblock In Hurd, J. \& Melham, T., editors, {\em TPHOLs}, LNCS 3603, pp.\
  245--260. Springer.

\bibitem[\protect\citeauthoryear{O'Connor}{O'Connor}{2009}]{oconnor-phd}
O'Connor, R. S.~S. (2009).
\newblock {\em Incompleteness \& Completeness: Formalizing Logic and Analysis
  in Type Theory}.
\newblock Ph.\ D. thesis, Radboud University Nijmegen.

\bibitem[\protect\citeauthoryear{Paulson}{Paulson}{2015}]{paulson-incompl-ar}
Paulson, L.~C. (2015).
\newblock A mechanised proof of {G\"odel's} incompleteness theorems using
  {Nominal Isabelle}.
\newblock {\em Journal of Automated Reasoning\/}~{\bf 55\/}(1), 1--37.
\newblock Online at
  \url{http://link.springer.com/article/10.1007%2Fs10817-015-9322-8}.

\bibitem[\protect\citeauthoryear{Shankar}{Shankar}{1986}]{shankar-phd}
Shankar, N. (1986).
\newblock {\em Proof-checking Metamathematics}.
\newblock Ph.\ D. thesis, University of Texas at Austin.

\bibitem[\protect\citeauthoryear{Shankar}{Shankar}{1994}]{shankar94}
Shankar, N. (1994).
\newblock {\em Metamathematics, Machines, and G{\"o}del's Proof}.
\newblock Cambridge University Press.

\bibitem[\protect\citeauthoryear{Shankar}{Shankar}{2013}]{shankar-deBruijn}
Shankar, N. (2013).
\newblock {Shankar, Boyer, Church-Rosser} and de {Bruijn} indices.
\newblock E-mail.

\bibitem[\protect\citeauthoryear{{\'S}wierczkowski}{{\'S}wierczkowski}{2003}]{swierczkowski-finite}
{\'S}wierczkowski, S. (2003).
\newblock Finite sets and {G{\"o}del's} incompleteness theorems.
\newblock {\em Dis\-ser\-ta\-tiones Math\-e\-ma\-ti\-cae\/}~{\bf 422}, 1--58.
\newblock \url{http://journals.impan.gov.pl/dm/Inf/422-0-1.html}.

\bibitem[\protect\citeauthoryear{Urban, \& Kaliszyk}{Urban \&
  Kaliszyk}{2012}]{urban-general}
Urban, C., \& Kaliszyk, C. (2012).
\newblock General bindings and alpha-equivalence in {Nominal Isabelle}.
\newblock {\em Logical Methods in Computer Science\/}~{\bf 8\/}(2:14), 1--35.

\bibitem[\protect\citeauthoryear{Wenzel}{Wenzel}{2007}]{wenzel-isabelle/isar}
Wenzel, M. (2007).
\newblock {Isabelle/Isar} --- a generic framework for human-readable proof
  documents.
\newblock {\em Studies in Logic, Grammar, and Rhetoric\/}~{\bf 10\/}(23),
  277--297.
\newblock From Insight to Proof --- Festschrift in Honour of Andrzej Trybulec.

\end{thebibliography}
\vspace*{10pt}

\address{COMPUTER LABORATORY\\
\hspace*{9pt}UNIVERSITY OF CAMBRIDGE\\
\hspace*{18pt}CAMBRIDGE, CB3 0FD, UK\\
{\it E-mail}: lp15@cam.ac.uk}
\end{document}